\documentclass[12pt]{article}

\newtheorem{defi}{Definition}

\newtheorem{tm}{Theorem}

\newtheorem{kor}{Corollary}
\newtheorem{rem}{Remark}

\usepackage{latexsym}
\usepackage{float}

\begin{document}

\begin{titlepage}


\begin{center}
{\Large\bf  Strongly regular graphs from orthogonal groups $O^+(6,2)$ and $O^-(6,2)$} 
\end{center}

\begin{center}
        Dean Crnkovi\'c (deanc@math.uniri.hr)\\
		Sanja Rukavina (sanjar@math.uniri.hr)\\
		and \\
		Andrea \v Svob (asvob@math.uniri.hr)\\[3pt]
		{\it\small Department of Mathematics} \\
		{\it\small University of Rijeka} \\
		{\it\small Radmile Matej\v ci\'c 2, 51000 Rijeka, Croatia}\\
\end{center}


\begin{abstract}
In this paper we construct all strongly regular graphs, with at most 600 vertices, admitting a transitive action of the orthogonal group $O^+(6,2)$ or $O^-(6,2)$.
Consequently, we prove the existence of strongly regular graphs with parameters (216,40,4,8) and (540,187,58,68).
We also construct a strongly regular graph with parameters (540,224,88,96) that was to the best of our knowledge previously unknown.
Further, we show that under certain conditions an orbit matrix $M$ of a strongly regular graph $\Gamma$ can be used to define a new strongly regular graph $\widetilde{\Gamma}$,
where the vertices of the graph $\widetilde{\Gamma}$ correspond to the orbits of $\Gamma$ (the rows of $M$). We show that some of the obtained graphs are related to each other
in a way that one can be constructed from an orbit matrix of the other.
\end{abstract}

\bigskip

{\bf AMS classification numbers:} 05E30, 05E18.


\end{titlepage}

\section{Introduction}

We assume that the reader is familiar with the basic facts of group theory and theory of strongly regular graphs.
We refer the reader to \cite{atlas, r} for relevant background reading in group theory, and to 
background reading in theory of strongly regular graphs we refer the reader to \cite{bjl, crc-srg, tonchev-book}.

The orthogonal groups $O^+(6,2)$ and $O^-(6,2)$ are simple groups of orders 20160 and 25920, respectively.
It is well known (see \cite{atlas}) that $O^+(6,2)=A_8=L(4,2)$ and $O^-(6,2)=O(5,3)=U(4,2)=S(4,3)$.

A graph is regular if all the vertices have the same degree; a regular graph is strongly regular of type $(v,k, \lambda , \mu )$ if it has $v$ vertices, degree $k$,
and if any two adjacent vertices are together adjacent to $\lambda$ vertices, while any two non-adjacent vertices are together adjacent to $\mu$ vertices.
A strongly regular graph  of type $(v,k, \lambda , \mu )$ is usually denoted by SRG$(v, k, \lambda, \mu)$. 
In this paper we classify strongly regular graphs with at most 600 vertices admitting a transitive action of the orthogonal group $O^+(6,2)$ or $O^-(6,2)$, 
using the method outlined in \cite{cms}. Among others, we construct strongly regular graphs with parameters (216,40,4,8) and (540,187,58,68), which proves the existence of strongly regular graphs
with these parameters. Further, one of the constructed SRGs with parameters (540,224,88,96) is new to our best knowledge.
Moreover, we introduce a method of constructing new strongly regular graphs from orbits, {\it i.e.} orbit matrices, of known strongly regular graphs.
In this paper we apply this method to the strongly regular graphs obtained from $O^+(6,2)$ and $O^-(6,2)$. 

\bigskip

For finding the graphs and computing their full automorphism groups, we used Magma \cite{magma} and GAP \cite{GAP2016, soicher}.

\section{SRGs from groups} \label{SRG_groups}

Using Theorem \ref{main} and Corollary \ref{main-graph} from \cite{cms}, we construct strongly regular graphs from orthogonal groups $O^+(6,2)$ and $O^-(6,2)$.

\begin{tm} \label{main}
Let $G$ be a finite permutation group acting transitively on the sets $\Omega_1$ and $\Omega_2$
of size $m$ and $n$, respectively.
Let $\alpha \in \Omega_1$ and $\Delta_2 =  \bigcup_{i=1}^s \delta_i G_{\alpha}$, where
$\delta_1,...,\delta_s \in \Omega_2$ are representatives of distinct $G_\alpha$-orbits.
If $\Delta_2 \neq \Omega_2$ and
$${\mathcal{B}}=\{ \Delta_2 g : g \in G \},$$
then ${\mathcal{D}}(G,\alpha,\delta_1,...,\delta_s)=(\Omega_2,{\mathcal{B}})$ is a
$1-(n, | \Delta_2 |, \frac{|G_{\alpha}|}{|G_{\Delta_2}|}\sum_{i=1}^{s} | \alpha G_{\delta_i} |)$ design with $\frac{m\cdot |G_{\alpha}|}{|G_{\Delta_2}|}$ blocks.
The group $H\cong G/{\bigcap_{x\in \Omega_2}G_x}$ acts as an automorphism group on $(\Omega_2,{\mathcal{B}})$, 
transitively on points and blocks of the design.

If $\Delta_2=\Omega_2$ then the set $\mathcal{B}$ consists of one block, and  ${\mathcal{D}}(G,\alpha,\delta_1,...,\delta_s)$ is
a design with parameters $1-(n,n,1)$.
\end{tm}

If a group $G$ acts transitively on $\Omega$, $\alpha \in \Omega$, and $\Delta$ is
an orbit of $G_{\alpha}$, then 
$\Delta' = \{ \alpha g \ | \ g \in G,\ \alpha {g^{-1}} \in \Delta \}$
is also an orbit of $G_{\alpha}$. $\Delta'$ is called the orbit of $G_{\alpha}$
paired with $\Delta$. It is obvious that $\Delta'' = \Delta$ and $| \Delta' | = | \Delta |$.
If $\Delta' = \Delta$, then $\Delta$ is said to be self-paired.

\begin{kor} \label{main-graph}
If $\Omega_1=\Omega_2$ and $\Delta_2$ is a union of self-paired and mutually paired orbits of $G_{\alpha}$, then the design 
${\mathcal{D}}(G,\alpha,\delta_1,...,\delta_s)$ is a symmetric self-dual design and the incidence matrix of that design is the 
adjacency matrix of a $|\Delta_2|-$regular graph. 
\end{kor}

The method of constructing designs and regular graphs described in Theorem \ref{main} is a generalization of results presented in \cite{cm, km, km1}. 
Using Corollary \ref{main-graph}, one can construct all regular graphs admitting a transitive action of the group $G$, but we will be interested only in those 
regular graphs that are strongly regular.

\subsection{SRGs from $O^+(6,2)$}

The alternating group $A_{8}$ is a simple group of order 20160, and up to conjugation it has 137 subgroups. 
There are 51 subgroups of the group $A_{8}$ up to index 600. In Table \ref{tb:subgrpsA8} we give the list of all the subgroups which lead to the construction of strongly regular graphs.
In the fifth column we give the rank of the permutation representation, {\it i.e.} the number of orbits of $H \le A_{8}$ acting on its cosets.
The sixth column indicates if the permutation representation is primitive or imprimitive.

\begin{table}[H]
\begin{center} \begin{footnotesize}
\begin{tabular}{|c|c|r|r|r|c|}
\hline
Subgroup& Structure& Order & Index & Rank & Primitive\\
\hline
\hline

$H^1_1$ & $S_6$ & 720 & 28 & 3 & yes \\
$H^1_2$ & $(A_4\times A_4):E_4$ & 576 & 35 & 3 & yes \\
$H^1_3$ & $E_4:(Z_2\times S_4)$ & 192 & 105 & 7 & no \\
$H^1_4$ & $(E_8:Z_7):Z_3$ & 168 & 120 & 4 & no \\
$H^1_5$ & $PSL(3,2)$ & 168 & 120 & 5 & no \\
$H^1_6$ & $E_9:D_8$ & 72 & 280 & 11 & no \\

\hline
\hline
\end{tabular}\end{footnotesize} 
\caption{\footnotesize Subgroups of the group $A_{8}=O^+(6,2)$}\label{tb:subgrpsA8}
\end{center}
\end{table}

Using the method described in Theorem \ref{main} and Corollary \ref{main-graph}, we obtained all the regular graphs on which the alternating group $A_{8}$ acts transitively and with at most 600 vertices. 
Using the computer search we obtained strongly regular graphs on 28, 35, 105, 120 or 280 vertices. Finally, we determined the full automorphism groups of the constructed SRGs. 

\begin{tm} \label{srg-A8}
Up to isomorphism there are exactly $6$ strongly regular graphs with at most $600$ vertices, admitting a transitive action of the group $A_{8}$.
These strongly regular graphs have parameters $(28,12,6,4)$, $(35,16,6,8)$, $(105,26,13,4)$, $(120,56,28,24)$ and $(280,117,44,52)$.
Details about the obtained strongly regular graphs are given in Table \ref{tb:srgA8}. 
\end{tm}

\begin{table}[H]
\begin{center} \begin{footnotesize}
\begin{tabular}{|c|c|c|}
\hline
Graph $\Gamma$ & Parameters  & $Aut (\Gamma) $  \\
\hline
\hline
$\Gamma^{1}_{1}=\Gamma(A_{8},H^1_{1})$ & (28,12,6,4) &  $S_8$\\ 
$\Gamma^{1}_{2}=\Gamma(A_{8},H^1_{2})$ & (35,16,6,8) &  $S_8$\\ 
$\Gamma^{1}_{3}=\Gamma(A_{8},H^1_{3})$ & (105,26,13,4) &  $S_{15}$\\ 
$\Gamma^{1}_{4}=\Gamma(A_{8},H^1_{4})$ & (120,56,28,24) &  $O^{+}(8,2):Z_2$\\ 
$\Gamma^{1}_{5}=\Gamma(A_{8},H^1_{5})$ & (120,56,28,24) &  $E_{64}:A_{8}$\\ 
$\Gamma^{1}_{6}=\Gamma(A_{8},H^1_{6})$ & (280,117,44,52) &  $S_9$\\ 

\hline
\hline
\end{tabular}\end{footnotesize} 
\caption{\footnotesize SRGs constructed from the group $A_{8}=O^+(6,2)$}\label{tb:srgA8}
\end{center}
\end{table}

\begin{rem} \label{IsoSrgs-A8}
The graphs $\Gamma^{1}_{1}$ and $\Gamma^{1}_{3}$ are the triangular graphs $T(8)$ and $T(15)$, respectively.
Moreover, $\Gamma^{1}_{3}$ is the unique strongly regular graph with parameters $(105,26,13,4)$.
Further, $\Gamma^{1}_{1}$ and $\Gamma^{1}_{2}$ are rank $3$ graphs, and $\Gamma^{1}_{4}$ is the complementary graph of the polar graph $NO^+(8,2)$.
$\Gamma^{1}_{5}$ is isomorphic to a strongly regular graph described in \cite{pg-mathon-street}, 
so its complement can be constructed from a partial geometry {\rm pg}$(7,8,4)$.
The graph $\Gamma^{1}_{6}$ is isomorphic to the {\rm SRG}$(280,117,44,52)$ constructed by R. Mathon and A. Rosa in \cite{m-r-srg}.
\end{rem}

\subsection{SRGs from $O^-(6,2)$}

The unitary group $U(4,2)$ is the simple group of order 25920, and up to conjugation it has 116 subgroups. 
There are 40 subgroups of the group $U(4,2)$ up to index 600. 
In Table \ref{tb:subgrpsU42} we give the list of all the subgroups which lead to the construction of strongly regular graphs.

\begin{table}[H]
\begin{center} \begin{footnotesize}
\begin{tabular}{|c|c|r|r|r|c|}
\hline
Subgroup& Structure& Order & Index & Rank & Primitive\\
\hline
\hline
$H^2_{1}$ & $E_{16}:A_{5}$ & 960 & 27 & 3 & yes \\
$H^2_{2}$ & $S_6$& 720 & 36 & 3 & yes \\
$H^2_{3}$ &$(E_9:Z_3):SL(2,3)$ & 648 & 40 & 3 & yes \\
$H^2_{4}$ &$E_{27}:S_4$ &648 & 40 & 3 & yes \\
$H^2_{5}$ &$(E_8.Z_{12}):Z_6$ & 576 & 45 & 3 & yes \\
$H^2_{6}$ &$(E_9:Z_3):Z_8$ & 216 & 120 & 7 & no \\
$H^2_{7}$ & $Ex_{32}^{+}:S_3$ & 192 & 135 & 9 & no \\
$H^2_{8}$ &$S_5$ & 120 & 216 & 10 & no \\
$H^2_{9}$ &$SL(2,3):Z_2$ & 48 & 540 & 25 & no \\
\hline
\hline
\end{tabular}\end{footnotesize} 
\caption{\footnotesize Subgroups of the group $U(4,2)=O^-(6,2)$}\label{tb:subgrpsU42}
\end{center}
\end{table}

Using the method described in Theorem \ref{main} and Corollary \ref{main-graph} we obtained all the regular graphs on which the unitary group $U(4,2)$ acts transitively and with at most 600 vertices.
Using the computer search we obtained strongly regular graphs on 27, 36, 40, 45, 120, 135, 216 or 540 vertices. Finally, we determined the full automorphism groups of the constructed SRGs.

\begin{tm} \label{srg-U(4,2)}
Up to isomorphism there are exactly $12$ strongly regular graphs with at most $600$ vertices, admitting a transitive action of the group $U(4,2)$.
These strongly regular graphs have parameters $(27,10,1,5)$, $(36,15,6,6)$, $(40,12,2,4)$, $(45,12,3,3)$, $(120,56,28,24)$, $(135,64,28,32)$, $(216,40,4,8)$, $(540,187,58,68)$ and $(540,224,88,96)$. 
Details about the obtained strongly regular graphs are given in Table \ref{tb:srgU42}. 
\end{tm}


\begin{table}[H]
\begin{center} \begin{footnotesize}
\begin{tabular}{|c|c|c|}
\hline
Graph $\Gamma$ & Parameters  & $Aut (\Gamma) $  \\
\hline
\hline
$\Gamma^{2}_{1}=\Gamma(U(4,2),H^2_{1})$ & (27,10,1,5) &  $U(4,2):Z_2$\\ 
$\Gamma^{2}_{2}=\Gamma(U(4,2),H^2_{2})$ & (36,15,6,6) &  $U(4,2):Z_2$\\ 
$\Gamma^{2}_{3}=\Gamma(U(4,2),H^2_{3})$ & (40,12,2,4) &  $U(4,2):Z_2$\\ 
$\Gamma^{2}_{4}=\Gamma(U(4,2),H^2_{4})$ & (40,12,2,4) &  $U(4,2):Z_2$\\ 
$\Gamma^{2}_{5}=\Gamma(U(4,2),H^2_{5})$ & (45,12,3,3) &  $U(4,2):Z_2$\\ 
$\Gamma^{2}_{6}=\Gamma(U(4,2),H^2_{6})$ & (120,56,28,24) &  $O^{+}(8,2):Z_2$\\ 
$\Gamma^{2}_{7}=\Gamma(U(4,2),H^2_{7})$ & (135,64,28,32) &  $O^{+}(8,2):Z_2$\\ 
$\Gamma^{2}_{8}=\Gamma(U(4,2),H^2_{8})$ & (216,40,4,8) &  $U(4,2):Z_2$\\ 
$\Gamma^{2}_{9}=\Gamma(U(4,2),H^2_{9})$ & (540,187,58,68) &  $Z_2\times(U(4,2):Z_2)$\\ 
$\Gamma^{2}_{10}=\Gamma(U(4,2),H^2_{9})$ & (540,187,58,68) &  $Z_2\times U(4,2)$\\ 
$\Gamma^{2}_{11}=\Gamma(U(4,2),H^2_{9})$ & (540,224,88,96) &  $Z_2\times U(4,2)$\\ 
$\Gamma^{2}_{12}=\Gamma(U(4,2),H^2_{9})$ & (540,224,88,96) &  $U(4,3):D_8$\\ 
\hline
\hline
\end{tabular} \end{footnotesize}
\caption{\footnotesize SRGs constructed from the group $U(4,2)=O^-(6,2)$}\label{tb:srgU42}
\end{center}
\end{table}

According to \cite{crc-srg, aeb}, the constructed strongly regular graph with parameters (216,40,4,8) and two strongly regular graphs with parameters (540,187,58,68) 
are the first known examples of strongly regular graphs with these parameters. Moreover, the constructed SRG(216,40,4,8) is the first known strongly regular graph on 216 vertices.
Further, the graph $\Gamma^{2}_{11}$ was to the best of our knowledge previously unknown.

\begin{rem} \label{IsoSrgs}
The strongly regular graphs $\Gamma^{1}_{4}$ and $\Gamma^{2}_{6}$ with parameters $(120,56,28,24)$ obtained from the groups $A_8$ and $U(4,2)$, respectively, are isomorphic.

The SRGs $(27,10,1,5)$, $(36,15,6,6)$, $(40,12,2,4)$ and $(45,12,3,3)$ are completely classified (see \cite{crc-srg, mckay-spence, spencesrg40, spencesrg45}).
The graphs $\Gamma^2_{1}$, $\Gamma^2_{2}$, $\Gamma^2_{3}$, $\Gamma^2_{4}$ and $\Gamma^2_{5}$ are rank $3$ graphs, where $\Gamma^{2}_{1}$ is the unique strongly regular graph on $27$ vertices.
Note that $\Gamma^{2}_{3}$ and $\Gamma^{2}_{4}$ are point graphs of generalized quadrangles {\rm GQ}$(3,3)$ (see \cite{h}), 
and $\Gamma^{2}_{4}$ corresponds to the point-hyperplane design in the projective geometry {\rm PG}$(3,3)$.
The graph $\Gamma^{2}_{5}$ is the only vertex-transitive strongly regular graph with parameters $(45,12,3,3)$.  
Further, $\Gamma^{2}_{7}$ is the complementary graph of the polar graph $O^+(8,2)$, and $\Gamma^{2}_{12}$ is the polar graph $NU(4,3)$.
\end{rem}

\begin{rem} \label{pg}
Strongly regular graphs can be constructed as point graphs of partial geometries (see \cite{crc-thas}). In particular, 
the existence of a partial geometry {\rm pg}$(11,16,4)$ would imply the existence of a {\rm SRG}$(540,187,58,68)$, but there is no known example of a partial geometry with these parameters.
If a strongly regular graph $\Gamma$ with parameters $(540,187,58,68)$ is obtained from a partial geometry {\rm pg}$(11,16,4)$, then a line in a {\rm pg}$(11,16,4)$ corresponds to a clique of size $12$ in $\Gamma$.
Since a {\rm pg}$(11,16,4)$ contains exactly $765$ lines, and $\Gamma^{2}_{10}$ has only $315$ cliques of size $12$, the graph $\Gamma^{2}_{10}$ cannot be obtained as the point graph of a {\rm pg}$(11,16,4)$.
The graph $\Gamma^{2}_{9}$ has exactly $1395$ cliques of size $12$, but there is no subset of $765$ cliques that can correspond to the lines of a {\rm pg}$(11,16,4)$.
Hence, the constructed strongly regular graphs with parameters $(540,187,58,68)$ cannot be obtained as point graphs of a {\rm pg}$(11,16,4)$, 
and the existence of a {\rm pg}$(11,16,4)$ remains undecided. In a similar way we conclude that the graph $\Gamma^{2}_{11}$ cannot be obtained as the point graph of a {\rm pg}$(14,15,6)$,
hence the existence of a {\rm pg}$(14,15,6)$ is still undetermined.
\end{rem}

\section{SRGs from orbit matrices}

M. Behbahani and C. Lam \cite{mb-lam} introduced the concept of orbit matrices of strongly regular graphs.
In \cite{dc-mm-bgr-sr}, the authors presented further properties of orbit matrices of strongly regular graphs.

Let $\Gamma$ be a SRG$(v, k, \lambda, \mu)$ and $A$ be its adjacency matrix. Suppose that an automorphism group $G$ of $\Gamma$ partitions the set of vertices $V$
into $t$ orbits $O_1, \ldots , O_t$, with sizes $n_1, \ldots , n_t$, respectively. 
The orbits divide the matrix $A$ into submatrices $[A_{ij}]$, where $A_{ij}$ is the adjacency matrix of vertices in $O_i$ versus those in $O_j$. 
We define the matrix $C = [c_{ij}]$, such that $c_{ij}$ is the column sum of $A_{ij}$. 
The matrix $C$ is the column orbit matrix of the graph $\Gamma$ with respect to the group $G$.
The entries of the matrix $C$ satisfy the following equations (see \cite{mb-lam, dc-mm-bgr-sr}):

\begin{equation}
\sum_{i=1}^t c_{ij}=\sum_{j=1}^t \frac{n_j}{n_i} c_{ij}=k,  \label{s3}
\end{equation}
\begin{equation}
\sum_{s=1}^t  \frac{n_s}{n_j} c_{is} c_{js}= \delta_{ij} (k- \mu ) + \mu n_i + ( \lambda - \mu) c_{ij}. \label{s4}
\end{equation}

While constructing strongly regular graphs with presumed automorphism group,
each matrix with the properties of a column orbit matrix, {\it i.e.} each matrix that satisfies equations (\ref{s3}) and (\ref{s4}),
is called a column orbit matrix for parameters $(v, k, \lambda, \mu)$ and orbith length distribution $(n_1, \ldots , n_t)$ (see \cite{mb-lam, dc-mm-bgr-sr}).
Hence, we have the following definition (see \cite[Definition~1]{dc-mm-bgr-sr}).

\begin{defi} \label{def-orb-mat}
A $(t \times t)$-matrix $C = [c_{ij}]$ with entries satisfying equations (\ref{s3}) and (\ref{s4})
is called a column orbit matrix for a strongly regular graph with parameters $(v,k, \lambda, \mu)$ and orbit lengths distribution $(n_1, \ldots ,n_t)$.
\end{defi} 

In the following theorem we show that, under certain conditions, a $(t \times t)$ orbit matrix of a SRG$(v,k, \lambda, \mu)$
can be used for a construction of a strongly regular graph on $t$ vertices.

\begin{tm} \label{om-construction}
Let $C= [c_{ij}]$ be a $(t \times t)$ column orbit matrix for a strongly regular graph $\Gamma$ with parameters $(v,k, \lambda, \mu)$ and orbit lengths distribution $(n_1, \ldots ,n_t)$,
$n_1= \ldots =n_t=n$, with constant diagonal. Further, let the off-diagonal entries of $C$ have exactly two values, {\it i.e.} $c_{ij} \in \{ x, y \}$, $x \neq y$, $1 \le i,j \le t$, $i \neq j$.
Replacing every $x$ with $1$ and every $y$ with $0$ in $C$, one obtains the adjacency matrix of a strongly regular graph $\widetilde{\Gamma}$ on $t$ vertices.  
\end{tm}
{\bf Proof} \enspace Since all the orbits have the same length, $C$ is a symmetric matrix.
Let all the diagonal elements be equal to $d$. Further, let us denote by $a$ the number of appearances of $x$ in the $i^{th}$ row.
Then the number of appearances of $y$ in the $i^{th}$ row is equal to $t-a-1$, and
$$ax+(t-a-1)y=k-d.$$
Obviously, $a$ does not depend on the choice of the selected row of $C$.

Let us take under consideration the $i^{th}$ and $j^{th}$ rows of $C$, and assume that $c_{ij}=x$. 
Denote by $b$ the number of columns in which the entry $x$ appear in both rows, {\it i.e.} $b=| \{ s \, | \, c_{is}=c_{js}=x  \} |$.
Then
$$\sum_{s=1}^t  c_{is}c_{js}=2dx+bx^2+2(a-b-1)xy+(t-2a+b)y^2,$$
which means that
$$2dx+bx^2+2(a-b-1)xy+(t-2a+b)y^2=\mu n + ( \lambda - \mu) x.$$
Since $b$ depends only on $x,y,d,a,n, \lambda$ and $\mu$, it does not depend on the choice of the rows $i$ and $j$, it depends only on the fact that $c_{ij}=x$.

Similarly, if $c_{ij}=y$ and we denote by $c$ the number of columns in which the entry $x$ appear in both rows, {\it i.e.} $c=| \{ s \, | \, c_{is}=c_{js}=x  \} |$,
then
$$\sum_{s=1}^t  c_{is}c_{js}=2dy+cx^2+2(a-c)xy+(t-2a+c-2)y^2.$$
Obviously, $c$ depends only on the value of $c_{ij}$, {\it i.e.} the fact that $c_{ij}=y$.
Hence, we proved that by replacing in $C$ every $x$ with $1$ and every $y$ with $0$ one obtains the adjacency matrix of a strongly regular graph 
with parameters $(t,a,b,c)$. 
$\quad \Box$

\subsection{Construction of SRGs from orbit matrices}

Applying the method given in Theorem \ref{om-construction}, we constructed strongly regular graphs from the orbit matrices of the graphs obtained in Section \ref{SRG_groups}. 
The results are presented in Table \ref{table-srg-OM}, following the notation from Theorem \ref{om-construction}. 
The third column contains information about the group for which the orbit matrix is constructed.

\begin{table}[htpb!]
\begin{center} \begin{footnotesize}
\begin{tabular}{|c | c | c | c | c | c | c |c|}
 \hline 
Graph $\Gamma$  & Parameters & $H \leq Aut (\Gamma) $& n & x & y & Graph $\widetilde{\Gamma}$ & Parameters
 \cr \hline \hline
 $\Gamma^{2}_{4}$ & (40,12,2,4) & $E_{8}$ & 4 & 0 & 2 &  $\widetilde{\Gamma^{2}_{4}}$& (10,3,0,1) \\
 \hline
$\Gamma^{1}_{4} \cong \Gamma^{2}_{6}$&(120,56,28,24)& $Z_{3}$& 3 &0 & 2& $\widetilde{\Gamma^{1}_{4}} \cong \widetilde{\Gamma^{2}_{6}} \cong \Gamma^{2}_{3}$&(40,12,2,4)\\
  \hline
 $\Gamma^{2}_{7}$&(135,64,28,32) & $Z_{3}$ & 3 &0 & 2& $\widetilde{\Gamma^{2}_{7}} \cong \Gamma^{2}_{5}$ & (45,12,3,3)   \\
  \hline \hline
\end{tabular} \end{footnotesize}
\caption{\footnotesize Strongly regular graphs constructed from orbit matrices} \label{table-srg-OM}
\end{center} 
\end{table}

\begin{rem} \label{1120}
The group $A_8$ acts transitively on the strongly regular graph $\Gamma_{1}={\rm SRG}(1120,390,146,130)$, which can be obtained from the group $O^+(8,3)$ as a rank $3$ graph. 
The full automorphism group of $\Gamma_{1}$ is isomorphic to $O^+(8,3):D_{8}$, and up to conjugation it has eighteen subgroups isomorphic to $Z_{4}$. 
The orbit matrix for one of these subgroups satisfies conditions given in Theorem \ref{om-construction}. 
The method described in Theorem \ref{om-construction} yields a {\rm SRG}$(280,36,8,4)$ isomorphic to the $U(4,3)$ polar graph listed in \cite{aeb}, 
having $U(4,3):D_{8}$ as the full automorphism group. Note that the group $A_8$ does not act transitively on the obtained strongly regular $(280,36,8,4)$ graph. 
Up to conjugation the group $Aut (\Gamma_{1})$ has five subgroups isomorphic to the group $Z_{28}$, and none of them produce an orbit matrix satisfying the conditions given in Theorem \ref{om-construction}. 
One of these subgroups, having all orbits of length $28$, yields an orbit matrix whose off-diagonal entries have exactly three values. Replacing one of these entries with $1$ 
and the other two with $0$, as described in {\it T}able \ref{table-srg-OMrem}, one obtains a strongly regular graph with parameters $(40,12,2,4)$. 
This is an example how the modification of the method described in Theorem \ref{om-construction} can be used for constructing strongly regular graphs, 
in case when some conditions of the theorem are not satisfied.

\begin{table}[htpb!]
\begin{center} \begin{footnotesize}
\begin{tabular}{|c | c | c | c | c | c | c |c|}
 \hline 
Graph $\Gamma$  & Parameters & $H \leq Aut (\Gamma) $& n & x & y & Graph $\widetilde{\Gamma}$ & Parameters
 \cr \hline \hline
 $\Gamma_{1}$ & (1120,390,146,130) & $Z_{4}$ & 4 & 4 & 1 &  $\widetilde{\Gamma_{1,1}}$& (280,36,8,4) \\
 \hline
&& $Z_{28}$& 28 &10 & 7, 13& $\widetilde{\Gamma_{1,2}} \cong {\Gamma^{2}_{4}}$&(40,12,2,4)\\
 
  \hline \hline
\end{tabular} \end{footnotesize}
\caption{\footnotesize Strongly regular graphs constructed from orbit matrices of a {\rm SRG}$(1120,390,146,130)$} \label{table-srg-OMrem}
\end{center} 
\end{table}

\end{rem}

\begin{rem} \label{om-theoretical}
Orbit matrices that can be used for the construction of SRGs applying the method introduced in Theorem \ref{om-construction}, do not have to be induced by an action of an automorphism
group on a strongly regular graph. It is sufficient that they are orbit matrices in terms of Definition \ref{def-orb-mat}. 
For example, the matrix $C_1$ from {\rm \cite{dc-mm-bgr-sr}} is a column orbit matrix for parameters $(40,12,2,4)$ and orbit length distribution $(4,4,4,4,4,4,4,4,4,4)$,
but it is not obtained as an orbit matrix induced by an action of an automorphism group on a {\rm SRG}$(40,12,2,4)$. 
However, $C_1$ has all the entries on the diagonal equal to $0$ and all the off-diagonal elements are equal to $0$ or $2$, hence 
replacing every off-diagonal $0$ with $1$ and every $2$ with $0$, we get the adjacency matrix of the Petersen graph, {\it i.e.} the strongly regular graph with parameters $(10,3,0,1)$. 
\end{rem}

\vspace*{0.2cm}

\noindent {\bf Acknowledgement} \\
This work has been fully supported by {\rm C}roatian Science Foundation under the project 1637. 


\end{document}